\documentclass[a4paper, 10pt]{article}
\usepackage{arXiv_style}

\usepackage{makecell}

\usepackage[colorlinks=true,linkcolor=red, citecolor=blue, urlcolor=blue]{hyperref}%

\newtheorem{theorem}{Theorem}[section]
\newtheorem{lemma}[theorem]{Lemma}
\newtheorem{corollary}[theorem]{Corollary}

\newtheorem{assumption}[theorem]{Assumption}

\newtheorem{definition}[theorem]{Definition}

\usepackage{wrapfig}

\usepackage{graphicx}

\usepackage{algpseudocode}
\usepackage{algorithm}
\usepackage{algorithmicx}

\usepackage{ amssymb }

\usepackage{url}
\usepackage{booktabs}       
\usepackage{nicefrac}       
\usepackage{microtype}      
\usepackage{natbib}
\usepackage{xcolor}

\usepackage{mathtools}

\usepackage{latexsym}
\usepackage{epsfig}
\usepackage{graphicx} 
\usepackage{dsfont}
\usepackage{mathtools}

\usepackage{caption}
\usepackage{subcaption}

\usepackage{amsmath, amsfonts, amssymb}
\usepackage{bbm, fancyhdr}
\usepackage{bm}
\usepackage{tikz}
\usepackage{tkz-euclide}
\usepackage{chngcntr}
\usepackage{verbatim}
\usepackage{mathtools}

\newcommand{\defeq}{\vcentcolon=}
\def\la{\langle}
\def\ra{\rangle}

\def \Z {\mathcal Z}

\def \X {\mathcal X}
\def \Y {\mathcal Y}
\def  \R {\mathbb R}
\def \E {\mathbb E}
\def \Gap {{\rm Gap}}

\renewcommand{\cite}[1]{\citep{#1}}

\title{
Bregman Proximal Method for Efficient Communications under Similarity}

\author{
  \bf{Aleksandr Beznosikov}
 \institute{MIPT, ISP RAS, Innopolis University}
  \and 
 \bf{ Darina Dvinskikh}
   \institute{HSE University}
  \and
\bf{ Dmitry Bylinkin}
 \institute{MIPT, ISP RAS}
    \and
\bf{Andrei Semenov}
\institute{MIPT}
    \and
\bf{Alexander Gasnikov}
\institute{Innopolis University, MIPT, ISP RAS}
     }

\begin{document}
\maketitle
%
%

\begin{abstract}
 We propose a novel stochastic distributed method for both monotone and strongly monotone variational inequalities with Lipschitz operator and proper convex regularizers arising in various applications from game theory to adversarial training. By exploiting \textit{similarity}, our algorithm overcomes the communication bottleneck that is a major issue in distributed optimization. 
The proposed method enjoys optimal communication complexity. All the existing distributed algorithms achieving the lower bounds under similarity condition essentially utilize the Euclidean setup. In contrast to them, our method is built upon the Bregman proximal maps and it is compatible with
an arbitrary problem geometry.
Thereby the proposed method fills an existing gap in this area of research. 
Our theoretical results are confirmed by numerical experiments on a stochastic matrix game.
\end{abstract}

\section{Introduction}

Variational inequalities (VIs) with monotone operators provide a unified and natural framework to cover many optimization problems, including convex minimization and convex-concave saddle point problems (SPPs) \citep{harker1990finite}. The interest in VIs is due to their wide applicability in  
economics, equilibrium theory, game theory, optimal control and differential equations, see \citep{facchinei2003finite,bauschke10convex}  for an introduction. They also play an important role in modern machine learning, see \citep{goodfellow2014generative, omidshafiei2017deep, madry2017towards, daskalakis2017training}.
We are interested in the regularized VI problem formulated as follows.
\begin{equation}\label{eq:VI}
\mbox{Find }  z^*\in \Z    :\quad  \langle F(z^*), z-z^* \rangle + g(z) - g(z^*)\geq 0, \quad\forall z\in\Z, 
\end{equation}
where   $F: \Z \to \R^d$ is a monotone and Lipschitz operator, $\Z\subseteq \R^d$ is a closed convex set, and $g:\Z\to\R$ is a proper convex lower semicontinuous function. Solving modern problems inevitably involves the use of huge training datasets. This forces engineers to utilize distributed computational systems. The data is distributed over $m$ nodes/clients/machines/devices coordinated by the server. Formally, this means working with an operator of the form $F(z) \defeq \frac{1}{m}\sum_{i=1}^m F_i(z)$. However, such coordination can significantly slow down the learning process, especially for systems with large computational resources \citep{bekkerman2011scaling}. Therefore, one of our main challenges in building the algorithm is to overcome the communication bottleneck. 

\subsection{ Similarity}
To reduce the communication frequency, many techniques have been developed. Among them is \textit{statistical preconditioning}, which reduces communication complexity by exploiting  \textit{similarity} \citep{shamir2014communication,hendrikx2020statistically}, i.e., the information that operators  $F_i$ are similar to each other and to their average $F$.
To account for similarity we will rewrite the VI problem \eqref{eq:VI} as a problem with two terms.
\begin{align}\label{eq:VInew}
    \mbox{Find }  z^*\in \Z    :\quad  \langle Q(z^*) + F_1(z^*), z-z^* \rangle +g(z)-g(z^*)\geq 0, \quad\forall z\in \Z, 
\end{align}
where $Q(z) \defeq F(z) - F_1(z)$ and $F_1$ is accessed by the server, which is the most computationally powerful node. The essence of similarity approaches is to move most of the computation to the server, offloading the other nodes. This allows to significantly reduce the amount of communication. In the case of convex minimization, this has been more than extensively studied.

\paragraph{Convex minimization.} The problem \eqref{eq:VI}  captures optimality conditions for constrained convex optimization with operators $F_i(z)= \nabla f_i(z)$. In machine learning applications, similarity is quite common: local functions $f_{i}$ are the average losses $\ell(\cdot, \cdot)$ on the training dataset $D_{i} = \{(x_1, y_1)^{(i)}, \dots, (x_N, y_N)^{(i)}\}$ stored at machine $i$:
\begin{equation}\label{ml_app}
f_{i}(z) \defeq \frac{1}{N} \sum_{j=1}^N \ell(z, (x_j, y_j)^{(i)}), \quad i=1\dots, m.
\end{equation}
Let us reformulate the minimization of the average risk on the functions from (\ref{ml_app}) as a stochastic optimization problem:
\begin{equation}\label{eq:sample_aver}
\min\limits_{z\in \Z} \left[f(z) \defeq \frac{1}{m}  \sum_{i=1}^m f_{i}(z)\right].
\end{equation}
where $f_{i}(z) \defeq f(z, D_i)$. $D_i$ is the local dataset, which can be considered as a set of random variables. If $D_{i}$'s on different machines are i.i.d. samples from the same distribution, the local empirical losses $f_{i}$ are statistically similar to their average $f$. There are various ways to measure similarity, but the most natural and theoretically justified  approach is the Hessian similarity: for all $ i =1,\dots, m$ with high probability  \cite{hendrikx2020statistically}
\begin{equation}\label{eq:hess_sim}
     \|\nabla^2 f(z) - \nabla^2 f_{i}(z)\| \leq \delta.
\end{equation}
This condition is referred to as  $\delta$-similarity, or $\delta$-relatedness. It is shown that $\delta$ can be estimated \cite{hendrikx2020statistically}. For this purpose, the objective $f$ is assumed to be $L$-smooth. Then, if the loss $f(z)$ is quadratic in $z$, then $\delta \sim \nicefrac{L}{\sqrt{N}}$ (up to a log factor), for a non-quadratic loss $\delta \sim \nicefrac{\sqrt{d} L}{\sqrt{N}}$ (up to a log factor) under the condition that 
\eqref{eq:hess_sim} holds uniformly over a compact domain with high probability \cite{zhang2018communication}. A similar analysis can be done for saddle point problems.

\paragraph{SPPs.}
The problem  \eqref{eq:VI} with operators 
$F_i(z) \defeq \left[\nabla_x f_i(x,y),~ - \nabla_y f_i(x,y)\right]^\top$ also captures convex-concave SPPs of the form
 \begin{equation}\label{eq:SPP_empir}
\min_{x \in \X} \max_{y \in \Y} \left[f(x,y)\defeq \frac{1}{m}\sum_{i=1}^m f_{i} (x,y)\right].   
\end{equation}
In this case, we measure  similarity in terms of second derivatives: for all $i=1,\dots, m$ with high probability  
\begin{align}\label{eq:spp_delta}
    & \| \nabla^2_{xx} f(x,y) - \nabla^2_{xx} f_{i}(x,y)  \| \leq \delta, \notag \\
    & \|\nabla^2_{xy} f(x,y) - \nabla^2_{xy} f_{i}(x,y)\| \leq \delta, \\
    &\|\nabla^2_{yy} f(x,y) - \nabla^2_{yy} f_{i}(x,y)\| \leq \delta.\notag
\end{align}
 A significant surge  of interest in SPPs is due to modern applications in GANs training \cite{gidel2018variational}, reinforcement learning \cite{jin2020efficiently,omidshafiei2017deep,wai2018multi}, distributed control \cite{necoara2011parallel}, and optimal transport \cite{jambulapati2019direct}.

\subsection{Related works}
VIs have been studied for more than half a century, and yet they remain an active area of research. Algorithms  solving VIs with Lipschitz and monotone operators date back to \citet{korpelevich1976extragradient}, who proposed the Extra Gradient method (an analog of the gradient method). A generalization of Extra Gradient  is Mirror Prox proposed by \citet{nemirovski2004prox}. The algorithm replaces the Euclidean projection with a more complex proximal Bregman step to fit the problem geometry.

\paragraph{Distributed methods.} 
 A distributed version of Extra Gradient was proposed by \citet{beznosikov2020distributed}. 
For the distributed version of Mirror Prox, see \citep{rogozin2021decentralized}. The authors also provided the lower bounds for the communication complexity: $\nicefrac{L}{\varepsilon}$ for the monotone VIs, and $\nicefrac{L}{\mu}$ for the strongly monotone VIs (up to a logarithmic factor). Here $\varepsilon$ measures the non-optimality gap function,  $L$ is the parameter of Lipschitz continuity, $\mu$ is the parameter of strong monotonicity.

\paragraph{Distributed methods exploiting  similarity.} 
 A first method exploiting Hessian similarity to reduce the communication complexity \textsc{DANE} was proposed by \citet{shamir2014communication}. It was designed for convex minimization problems. Later $\Omega\left( \sqrt{\nicefrac{\delta}{\mu}}\log\nicefrac{1}{\varepsilon} \right)$ was established as the lower boundary of communication rounds under Hessian similarity \citep{arjevani2015communication}. For quadratic problems, an optimal method has been developed using a
modification of the DANE approach \cite{yuan2020convergence}. 
For an optimal (up to the log factor)  method in the general case, see \citep{tian2022acceleration}. For SPPs, optimal (up to a log factor) algorithms have been proposed by
 \citet{beznosikov2021distributed} together with 
 the lower bounds of communication complexity: $\nicefrac{\delta}{\varepsilon}$ for the convex-concave SPPs, and $\nicefrac{\delta}{\mu}\cdot\log\nicefrac{1}{\varepsilon}$ for the strongly convex-strongly concave SPPs. \citet{kovalev2022optimal} improved the results of \citet{tian2022acceleration} and  \citet{beznosikov2021distributed} by proposing optimal methods for convex minimization and convex-concave SPPs with optimal computational complexity, i.e. the authors eliminated non-optimal logarithm. In addition, there are a number of state-of-the-art methods that work with similarity, the complexity of which is improved by using communication compression and client sampling \citep{beznosikov2024similarity, beznosikov2022compression, lin2024stochastic}.

\subsection{Contribution}
All the existing distributed algorithms exploiting similarity for   convex minimization, SPPs or VIs  are substantially utilizes the Euclidean setup.
However, a large number of problems have non-Euclidean geometry, e.g., minimization on the probability simplex  arising  in machine learning \cite{nemirko2021machine},  statistics, chemistry, portfolio management \cite{chang2000heuristics},  optimal transport and Wasserstein barycenters \cite{agueh2011barycenters}. 
Indeed, the Euclidean distance is not  well suited for probability measures.  Generalizing Euclidean algorithms to non-Euclidean ones is   non-trivial and usually requires the development of new methods. Motivated by this gap between the Euclidean and non-Euclidean algorithms, we aim to design a novel method to tackle the problem geometry.
To the best of our knowledge, we present the first distributed method utilizing similarity with the non-Euclidean setup.
Technically our approach is based on the Bregman proximal maps, which work particularly well in constrained optimization.
Our contribution  can be summarized as follows:
\begin{itemize}
    \item We present  a \textbf{P}roximal \textbf{A}lgorithm \textbf{u}nder \textbf{S}imilarity (\textsc{PAUS}) that utilizes similarity  for VIs with a monotone Lipschitz operator and a convex composite $g$. It achieves optimal communication complexity of $\nicefrac{\delta}{\varepsilon}$, where $\varepsilon$ measures the non-optimality gap function, and $\delta$ is the parameter of similarity.
    \item Our analysis also shows a speedup of \textsc{PAUS} in the case of an operator, that is $\mu$-strongly monotone with respect to the Bregman divergence. Our method achieves optimal communication complexity of $\nicefrac{\delta}{\mu}\cdot\log\nicefrac{1}{\varepsilon}$ up to a logarithmic factor. The analysis is compatible with an arbitrary Bregman divergence.
    \item In both cases, we generalize our analysis by adding stochasticity to the method. \textsc{PAUS} achieves $\nicefrac{\delta}{\varepsilon}+\nicefrac{\sigma_*^2}{\varepsilon^2}$ for monotone VIs and $\nicefrac{\delta}{\mu}\cdot\log\nicefrac{1}{\varepsilon} + \nicefrac{\sigma_*^2}{\mu\varepsilon}$ for strongly monotone ones. In contrast to classical works that use the \textsc{SGD}-like approach to construct stochastic methods, we require boundedness of the variance of the stochastic oracle $\sigma_*^2$ only at the solution.
    \item  We confirm our theoretical results by numerical experiments.
\end{itemize}

\paragraph{Paper organization.}
The structure of the paper is as follows. 
Section \ref{sec:prel} introduces the preliminaries: the necessary definitions and assumptions. Section \ref{sec:alg} presents the stochastic algorithm for solving VIs with monotone and Lipschitz operators and convex composites. Section \ref{sec:strongly} extends the convergence theory to the case of strongly monotone VIs. Section \ref{sec:exper} presents the numerical experiments demonstrating the superiority of our novel algorithm on the two-player stochastic matrix game.

\section{Preliminaries}\label{sec:prel}
\paragraph{Notation.} For vectors, $\|z\|$ is a general norm on space $\Z$, and $\|s\|_{*}$ is its dual norm on the dual space $\Z^*$: $ \|s\|_{*} = \max_{z\in \Z} \{ \la s,x \ra : \|z\| = 1 \} $. For matrices,  $\|A\|$ is the matrix norm induced by vector norm $\|z\|$: $\|A\| = \sup_{z\in \Z}\{\|Az\| : \|z\| =1 \}.$

When discussing the communication efficiency, it is important to choose the right definition. One possible quality metric is the number of communication rounds. In addition to this, each vector exchange between client and server or the time spent on communication can be considered. This paper assumes the definition of communication complexity in the first mentioned sense. We measure how often the server communicates with the nodes without considering the number of vector exchanges per round.

Methods that utilize data similarity move the most computational complexity to the server. Therefore, in addition to communication complexity, it is worth measuring the complexity of local computations performed on the server.

Next, we provide all the definitions and assumptions necessary to build the convergence theory.

\begin{definition}\label{def:strongly}
    We say that $g:\Z \to \R$ is a $\mu$-strongly convex function with respect to $\|\cdot\|$, if
    \begin{align*}
        g(u)-g(v)\geq\langle h,u-v\rangle+\frac{\mu}{2}\|u-v\|^2,\quad\forall u,v \in Z,\quad h\in\partial g(v).
    \end{align*}
\end{definition}
If $\mu=0$, $g$ is called a convex function.

\begin{definition}
    We say that $w:\Z \to \R$ is a distance generating function (DGF), if $w$ is 1-strongly convex with respect to $\|\cdot\|$, i.e., for all $u, v \in \Z$: $w(v) \geq w(u) + \langle \nabla w(u), v - u\rangle +\frac{1}{2}\|v-u\|^2$. The corresponding Bregman divergence is 
    \[
    V(u,v) = w(u) - w(v) -\langle \nabla w(v), u-v\rangle, \quad \forall u, v \in \Z.
    \]
\end{definition}
The property of  distance generating function ensures $V(u,u)=0$ and $V(u,v)\geq \frac{1}{2}\|u-v\|^2$. 

Since we aim to use stochastic methods to reduce computation time, we introduce $F(\cdot,\xi)$ and $F_1(\cdot,\xi)$. These are the stochastic oracles of the corresponding operators. Consider $F(z,\xi)$. The random variable $\xi$ can be understood in different ways. For example, as a sample of the computational node to be communicated with: $F(z,\xi)=F_{\xi}(z)$. Another interesting case covered by our analysis is the imposition of additive noise on the operator belonging to each machine: $F(z,\xi)=\frac{1}{m}\sum_{i=1}^m\left( F_i(z)+\xi_i \right)=F(z)+\xi$. Informally, such a technique allows to maintain privacy protection of local data. This refers to Federated Learning, a very popular trend nowadays \citep{li2020federated}. We take into account that the objective operator $F(\cdot)$ has the form of a finite sum, where each $F_i(\cdot)$ can also be a finite sum. Thus, it is possible to speed up the server by using a stochastic approach: $F_1(z,\xi)=F_{1,\xi}(z)$.

\begin{assumption}[Stochastic oracle for $F_1(\cdot)$]\label{ass:stochastic_server}
    The stochastic oracle $F_1(\cdot,\xi)$ is unbiased and its variance is bounded at the solution:
    \begin{align*}
        \E_{\xi}[F_1(z,\xi)]=F_1(z),\quad \E_{\xi}[\|F_1(z^*,\xi)-F_1(z^*)\|_*^2]\leq\sigma_{1,*}^2,\quad\forall z\in \Z,
    \end{align*}
\end{assumption}

\begin{assumption}[Stochastic oracle for $F(\cdot)$]\label{ass:stochastic}
    The stochastic oracle $F(\cdot,\xi)$ is unbiased and there are two options of variance boundedness:
    \begin{align*}
        \E_{\xi}[F(z,\xi)]=F(z)\quad\forall x\in\R.
    \end{align*}
    \paragraph{(a)} The variance of $F(\cdot,\xi)$ is uniformly bounded:
    \begin{align*}
        \E_{\xi}[\|F(z,\xi)-F(z)\|_*^2]\leq\sigma^2,\quad \forall z\in \Z.
    \end{align*}
    \paragraph{(b)} The variance of $F(\cdot,\xi)$ is bounded at the solution:
    \begin{align*}
        \E_{\xi}[\|F(z^*,\xi)-F(z^*)\|_*^2]\leq\sigma_*^2,\quad \forall z\in \Z.
    \end{align*}
\end{assumption}
Here the uniform boundedness is only needed to show convergence by the dual gap function in the monotone case. For strongly monotone VIs, only boundedness at the solution is enough.

\begin{assumption}[Monotonicity]\label{ass:monotone}
     The operators $F(\cdot,\xi)$, $F_1(\cdot,\xi)$ are monotone, i.e. for all $u, v \in \Z$ and for every $\xi$:
\[
\langle F(u,\xi)-F(v,\xi), u - v \rangle \geq 0,
\]
\[
\langle F_1(u,\xi)-F_1(v,\xi), u - v \rangle \geq 0.
\]
\end{assumption}
A special case of monotone VIs are convex optimization problems and convex-concave SPPs.  

\begin{assumption}[Lipschitzness]\label{ass:Lipsch}
    The operator $F_1(\cdot,\xi)$ is $L_{F_1}$-Lipschitz continuous, i.e. for all $u, v \in \Z$ and for every $\xi$:
    \[
    \|F_1(u,\xi) - F_1(v,\xi)\|_* \leq L_{F_1}\|u - v\|.
    \]
\end{assumption}
For convex minimization, this means that server function $f_1(z,\xi)$ is $L_{F_1}$-smooth.

\begin{assumption}[$\delta$-similarity]\label{ass:delta}
    The operator $F_1(\cdot) - F(\cdot,\xi)$ is $\delta$-Lipschitz continuous, i.e., for all $u, v \in \Z$ :
\[ 
\| F_1(u) - F(u,\xi) - F_1(v) + F(v,\xi)  \|_* \leq \delta \|u - v\|.
\]
\end{assumption}
In the Euclidean case, this is a stochastic generalization of $\delta$-similarity for convex minimization and convex-concave SPPs. This is evident by considering \eqref{eq:hess_sim} and \eqref{eq:spp_delta} and noting that uniform bounding the norm of the Hessian by some constant $\delta$ entails $\delta$-smoothness of the function. Note that all assumptions are introduced on the stochastic oracles of the operators, which is a stronger case than if we had introduced them on the operators themselves. This is necessary in order to use a wider class of stochastic oracles whose variance is bounded only at the optimum. If the variance of the chosen stochastic oracle is uniformly bounded, then we can relax Assumptions \ref{ass:monotone} and \ref{ass:Lipsch}, but not Assumption \ref{ass:delta}.

\section{Monotone case}\label{sec:alg}
In this section, we present the algorithm which solves \eqref{eq:VInew}.

\subsection{Main algorithm}

Now we are in a position to provide \textsc{PAUS} for VIs, see Algorithm \ref{Alg:VI}.

\begin{algorithm}[ht!]
\caption{\textsc{PAUS}}
\label{Alg:VI}  
\begin{algorithmic}[1]
\Require  parameter of similarity $\delta$, stepsize $\gamma \leq  \nicefrac{1}{2\delta}$, parameter $\alpha\geq0$, number of iterations $K$, starting points $z^0=u^0\in \Z$
\For{$k= 0,1,2,\dots, K-1$}
\State \hspace{-0.3cm} Sample random variable $\xi^k$ on server
 \State \hspace{-0.3cm} Collect $F(z^k,\xi^k)=\frac{1}{m}\sum_{i=1}^mF_i(z^k,\xi_i^k)$ on server\label{compute_zk}
 \State \hspace{-0.3cm} Find $u^k$ as a solution to \label{line:SPP_prob} 
  \[ \langle \gamma  (   F_1(u^k) +F(z^k,\xi^k) -F_1(z^k))+ \nabla w (u^k)  - \nabla w (z^k), z - u^k \rangle + \gamma\left(g(z)-g(u^k)\right) \geq 0\] 
  \Statex \hspace{0.25cm} for all $z\in \Z$ by \textsc{SCMP} (Algorithm \ref{Alg:CMP}) procedure on server 
\State \hspace{-0.3cm} Collect $F(u^k,\xi^k)=\frac{1}{m}\sum_{i=1}^mF_i(u^k,\xi_i^k)$ on server\label{compute_uk}
\State \hspace{-0.3cm} Solve \label{eq:line5Paus}
\[ z^{k+1} = \arg\min_{z \in \Z} \left\{ \gamma \langle F(u^k,\xi^k) - F_1(u^k)  -F(z^k,\xi^k) + F_1(z^k), z \rangle+ (1+\alpha)V(z,u^k) \right\}\]
\Statex \hspace{0.25cm} on server

\EndFor
 \State \textbf{return} $ \widetilde{u}^K = \frac{1}{K}\sum_{k=0}^{K-1} u^k$ for monotone VIs and $z^K$ for strongly monotone ones
\end{algorithmic}
\end{algorithm}

At each iteration $k = 0,1,\dots, K-1$ of Algorithm \ref{Alg:VI}, the server initiates two rounds of communication to compute $F(z^k,\xi^k)$ at Line \ref{compute_zk} and $F(u^k,\xi^k)$ at Line \ref{compute_uk}. To solve the inner problem  encountered in Line \ref{line:SPP_prob} we provide a procedure we call  \textsc{Stochastic Composite MP (SCMP)} (see Algorithm \ref{Alg:CMP}). Importantly, the same random variable $\xi^k$ is used in the formulation of both subproblems. This is inspired by the work of \citet{mishchenko2020revisiting}, where this technique is used to not require uniform boundedness of the stochastic oracle.

We further propose a descent lemma for \textsc{PAUS}. 

\begin{lemma}\label{descent_lemma}
    Consider Assumption \ref{ass:delta}. Then the inequality
    \begin{align*}
  2\gamma \left[\langle F(u^k,\xi^k), u^k - z \rangle + g(u^k)-g(z)\right] \leq& 2V(z, z^k) - 2V(u^k, z^k) - 2(1+\alpha)V(z,z^{k+1}) \notag \\
  &- 2V(z^{k+1}, u^k)+2\alpha V(z,u^k) \notag+\gamma^2\delta^2\|u^k-z^k\|^2 \notag\\&+\|z^{k+1}-u^k\|^2.
    \end{align*}
    holds.
\end{lemma}
See the proof in Appendix \ref{app:descent_lemma}. The next theorem presents the convergence rate of \textsc{PAUS} by the following gap function: 
\begin{equation}\label{eq:gap}
    \Gap(u)  = \max_{z \in \Z} \left\{  \langle F(z),u-z \rangle +g(u)-g(z)\right\}.
\end{equation}
This function is the standard criterion for VIs. It corresponds to the standard optimality criteria in convex minimization and SPPs \cite{nemirovski2004prox,juditsky2011solving}.

\begin{theorem}\label{th:PAUS}
Consider assumptions of Lemma \ref{descent_lemma} with Assumptions \ref{ass:stochastic}(a) and \ref{ass:monotone}. Then after $K $ communication rounds, \textsc{PAUS} (Algorithm \ref{Alg:VI}),  run with $\alpha=0$, stepsize $\gamma \leq\nicefrac{1}{2\delta}$ and starting points $z^0,u^0\in \Z$,
outputs  $\tilde{u}^K$ such that
\begin{align*}
    \E\left[\Gap(\tilde{u}^K)\right]\leq\frac{D}{K\gamma}+ \frac{2\gamma}{3}\sigma^2,
\end{align*}
where $D=\sup_{z\in Z}\left\{V(z, z^0)\right\}$.
\end{theorem}
See the proof in Appendix \ref{app:lemma_PAUS}. The theorem guarantees convergence of the method to some neighborhood of the solution. The following corollary shows that if $\gamma$ is chosen correctly, convergence with arbitrary accuracy can be guaranteed.

\begin{corollary}\label{cor:monotone}
 Consider assumptions of Theorem \ref{th:PAUS}.
Let $\tilde{u}^K$ be the output of \textsc{PAUS} (Algorithm \ref{Alg:VI}), run with appropriate parameters and starting points $z^0,u^0\in \Z$ 
in
 \[ \mathcal{O}\left( \frac{D\delta}{\varepsilon}+\frac{D\sigma^2}{\varepsilon^2} \right) \text{communication rounds,}
  \]  
then  $\Gap(\tilde{u}^K)\leq \varepsilon$. 
\end{corollary}
See proof in Appendix \ref{proof:cor:monotone}. 

\subsection{Stochastic Composite MP}
Let us rewrite the problem encountered in Line \ref{line:SPP_prob}  of Algorithm \ref{Alg:VI} as finding $v^*$ such that for all $z\in \Z$:
\begin{align}\label{eq:prob-COMMP}
  \langle \gamma  (  F_1(v^*) + F(z^k,\xi^k) - F_1(z^k))+ \nabla w (v^*)  - \nabla w (z^k), z - v^* \rangle+\gamma(g(z)-g(v^*))\geq 0.  
\end{align}
Note that in Line \ref{line:SPP_prob} we have to solve the problem with a strongly monotone operator, i.e., the Bregman divergence can also be used as a convergence criterion.

\begin{algorithm}[ht!]
\caption{ to solve auxiliary problem in step \ref{line:SPP_prob}  of Algorithm \ref{Alg:VI} }
\label{Alg:CMP}  
\begin{algorithmic}[1]
\Procedure{\textsc{SCMP}}{$\gamma, L_{F_1},  z^k$}
\State Choose stepsize $\eta =\frac{1}{3\gamma L_{F_1}}$ 
\State Choose starting point $v^0 \in \Z$
\For{$t= 0,1,2,\dots, T-1$}
 \State \hspace{-0.5cm} Sample random variable $\xi^t$ on server
 \State \hspace{-0.5cm} Solve $v^{t+\frac{1}{2}} = \arg\min\limits_{v\in \Z} \{ \gamma\eta\langle H(v^t,\xi^t), v  \rangle 
  + \eta V(v, z^k) + V(v, v^t)+\gamma g(v)\}$ on server\label{eq:MP_line5}
\State \hspace{-0.5cm} Solve $v^{t+1} = \arg\min\limits_{v\in \Z} \{ \gamma\eta\langle H(v^{t+\frac{1}{2}},\xi^t), v  \rangle  + \eta V(v, z^k) + V(v, v^t) +\gamma g(v)\}$ on server\label{eq:MP_line6}
\EndFor
 \State \textbf{return} $v^{T}$
  \EndProcedure
\end{algorithmic}
\end{algorithm}

The next theorem presents the convergence guarantee for \textsc{SCMP} to solve the problem \eqref{eq:prob-COMMP}. For the sake of brevity of description, we denote $H(v,\xi)=\gamma  (  F_1(v,\xi) + F(z^k,\xi^k) - F_1(z^k))$.
\begin{theorem}\label{th:compositeMP}
   Let Assumptions \ref{ass:stochastic_server}, \ref{ass:monotone}, \ref{ass:Lipsch} and \ref{ass:delta} hold. Then \textsc{SCMP} (Algorithm \ref{Alg:CMP}) with $\gamma=\frac{1}{2\delta}$ produce the sequence $\{v^t\}$ such that 
   \begin{align*}
    \E\left[V(v^*, v^{t+1})\right] \leq \left(1-\frac{\eta}{2}\right)\E\left[V(v^*,v^t)\right] + 4\eta^2\sigma_{1,*}^2.
    \end{align*}
\end{theorem}
The proof of this theorem is given in Appendix 
\ref{app:irror_comp}. As in the case of Theorem \ref{th:PAUS}, fine tuning of the algorithm parameters is required to obtain convergence with arbitrary accuracy.
\begin{corollary}\label{cor:scmp}
   Consider assumptions of Theorem \ref{th:compositeMP}. Let $ v^{T}$  be the output of  \textsc{SCMP} procedure. Consider stepsize $\gamma=\nicefrac{1}{2\delta}$ and starting point $v^0$. Then Algorithm \ref{Alg:CMP} with appropriate choice of $\eta$ needs
   $$
   \mathcal{O}\left(\frac{L_{F_1}}{\delta}\log\frac{V(v^*,v^0)}{\varepsilon} + \frac{\sigma_{1,*}^2}{\varepsilon}\right)
   \text{ iterations} 
   $$ 
   to achieve $V(v^*, v^T)\leq \varepsilon$.
\end{corollary}
See proof in Appendix \ref{proof:cor:scmp}. The procedure supports the possibility of a stochastic solution. This can be useful when the server node has too much data and needs to speed up the computation. We use the "same sample" technique from \citep{mishchenko2020revisiting}. To solve the subproblem by the proposed method, we only need boundedness of the variance of the stochastic oracle $F_1(\cdot,\xi)$ at the solution.

\section{Extension to strongly monotone VIs}\label{sec:strongly}
In this section, we obtain a linear convergence rate of \textsc{PAUS} by the Bregman divergence, which is an appropriate convergence criterion for the strongly monotone operator.

\begin{assumption}\label{strong_monotonicity}
    Operator $F(\cdot,\xi)$ is $\mu$-strongly monotone with respect to \textsc{DGF} $w(\cdot)$, i.e. for all $u,v\in Z$ and for every random variable $\xi$:
    \begin{align*}
        \langle F(u,\xi)-F(v,\xi),u-v \rangle\geq \frac{\mu}{2}\left(V(u,v)+V(v,u)\right).
    \end{align*}
\end{assumption}
We assume strong monotonicity with respect to the corresponding geometry. This definition is not used for the first time, it is often found in the literature \citep{lu2018relatively, ablaev2022some, stonyakin2021some}.

The main difficulty in constructing a convergence theory for strongly monotone VIs is that the Bregman divergence is in general non-symmetric and there is no triangle inequality for it. To solve this problem, we "inflate" the coefficient in the subproblem (\ref{eq:line5Paus}) by some value $\alpha$ and use the strong monotonicity of the operator to eliminate the extra summands.

\begin{theorem}\label{th:strongly_monotone}
Consider assumptions of Lemma \ref{descent_lemma} with Assumption \ref{ass:stochastic}(b) and Assumption \ref{strong_monotonicity}.
Consider $\alpha=\nicefrac{\gamma\mu}{2}$, $\gamma \leq \nicefrac{1}{2\delta}$ and a starting point $z^0\in \Z$. Then the inequality
\begin{align*}
    \E\left[V(z^*,z^{k+1})\right]\leq\left(1-\frac{\gamma\mu}{4}\right)\E\left[V(z^*,z^k)\right]+\frac{2\gamma^2}{3}\sigma_*^2
\end{align*}
holds.
\end{theorem}
See the proof in Appendix \ref{app:strongly}. Let us repeat the proof of Corollary \ref{cor:scmp} with other constants and obtain
\begin{corollary}
 Consider assumptions of Theorem \ref{th:strongly_monotone}.
Let $z^K$ be the output of \textsc{PAUS} (Algorithm \ref{Alg:VI}), run  with an appropriate parameters and a starting point $z^0\in \Z$, 
in
 \[ \mathcal{O}\left( \frac{8\delta}{\mu}\log\frac{1}{\varepsilon}+\frac{8\sigma_*^2}{3\mu \varepsilon} \right) \text{ communication rounds,}
  \]  
then  $V(z^*,z^K)\leq \varepsilon$.    
\end{corollary}
Note that \textsc{PAUS} performs two rounds of communication per iteration. Thus, asymptotically the communication complexity by rounds is equal to the complexity by iterations.

In \textsc{PAUS} (Algorithm \ref{Alg:VI}), the main computational load is shifted to the server. At one iteration of the method, each device accesses the local oracle twice, while the server is forced to do it $\mathcal{O}\left(\nicefrac{L_{F_1}}{\delta}\cdot\log\nicefrac{V(v^*,v^0)}{\varepsilon} + \nicefrac{\sigma_{1,*}^2}{\varepsilon}\right)$ times. Nevertheless, the additional freedom to choose the Bregman divergence can make it easier to compute the proximal mapping or even produce a closed-form solution. See examples in Appendix \ref{closed_forms_app}.

\section{Numerical experiments}\label{sec:exper}


We  evaluate the effectiveness of  \textsc{PAUS}
by comparing it with other distributed algorithms with and without exploiting data similarity. As a first algorithm
we consider \textsc{decentralized Mirror Prox} from \citep{rogozin2021decentralized}. This algorithm is based on the Bregman proximal maps but does not exploit similarity.  
The second algorithm  is the algorithm from \citep{kovalev2022optimal} which we will refer to as \textsc{Extra Gradient}. The algorithm  
is designed under similarity but in the Euclidean setup.

\paragraph{Problem.}  We carry out numerical experiments for a stochastic matrix game 
\[  \min\limits_{x \in \Delta} \max\limits_{y \in \Delta}\left[ x^\top \E [A_{\xi}] y\right],\]
where $x, y$ are the mixed strategies of  two players, $\Delta$ is the probability simplex,  and $A_{\xi}$ is a stochastic payoff
matrix. A solution of this problem can be approximated by a solution of the following empirical problem
\begin{equation*}\label{eq:matrix_game}
  \min_{x \in \Delta} \max_{y \in \Delta} \left[\frac{1}{m} \sum_{i=1}^m x^\top  A_{i} y\right],
\end{equation*}
where $A_{1}, \dots, A_m$ are 
i.i.d samples  of stochastic matrix  $A_{\xi}$.
We study the convergence of algorithms in terms of
 the duality gap 
\[\max\limits_{y \in \Delta} \left[\frac{1}{m}\sum_{i=1}^m \left(\widetilde x^K\right)^\top  A_i y\right]  - \min\limits_{ x \in \Delta}  \left[\frac{1}{m}\sum_{i=1}^m  x^\top  A_i \widetilde y^K\right]   \leq \varepsilon,\]
 where $\widetilde u^K = (\widetilde x^K, \widetilde y^K)$ is the output of algorithms. This duality gap is the upper bound for the gap function for VIs from \eqref{eq:gap}.  

\begin{wrapfigure}[11]{r}{0.4\textwidth}
\vspace{-1.4cm}
\includegraphics[width=0.4\textwidth]{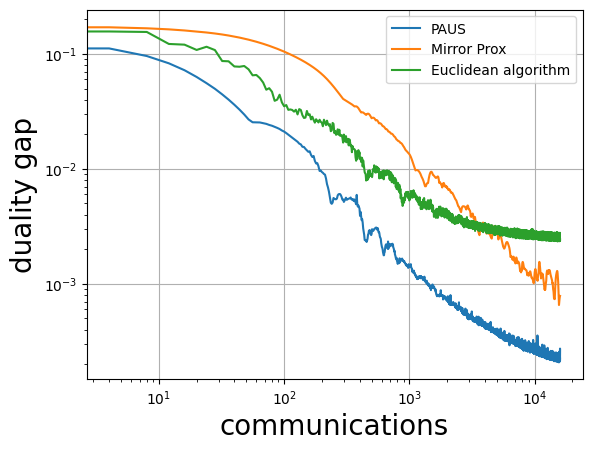}
\vspace{-0.8cm}
\caption{ Comparison of \mbox{state-of-the-art} methods  }
\label{fig:conver_comp}
\end{wrapfigure}

 \paragraph{Data.}
 We model entries of stochastic matrix $A_\xi$
of size $d\times d$ with $d = 25$ as 
$ [A_\xi]_{ij} = (1+\nu\xi)[C]_{ij},$
 where    $\nu\in(0,1)$, $\xi$ is a Rademacher random variable and  $C$ is a deterministic matrix given in Exercise 5.3.1 of \citep{ben2001lectures}. $\nu$ is used to adjust the similarity parameter $\delta$.
From this matrix we sample  $m = 10^4$ matrices $A_1,\dots, A_m$. These matrices are shared between  $5$ devices, each holds  local datasets  of size  $n =2\cdot 10^3$.  

Figure \ref{fig:conver_comp} demonstrates the superiority of \textsc{PAUS} in comparison with other distributed algorithms on the problem with $L\approx10^{-1}$ and $\delta\approx10^{-2}$, $\nu\approx10^{-3}$. The parameters of algorithms were estimated theoretically (see Appendix \ref{sec:additional_ex}) and then tuned to get faster convergence. 
All algorithms have approximate slope ratio $-1$ according to their theoretical bounds ($K \sim \varepsilon^{-1}$). A faster convergence of our algorithm in comparison with the Euclidean algorithm \cite{kovalev2022optimal}  is achieved due to better 
utilizing the constrained set, namely the probability simplex, for which the Euclidean distance is worse than the $\ell_1$-distance. Thus our method is able to better estimate the parameter $\delta$  in the proper norm and the distance to a solution (in terms of the Bregman divergence).
A slow convergence of Mirror Prox \cite{rogozin2021decentralized}
 is explained by ignoring data similarity.





\section{Conclusion}
In this paper, we introduced the novel distributed stochastic proximal algorithm using similarity for both monotone and strongly monotone variational inequalities: \textsc{PAUS}. It achieves optimal communication complexity:  
$\nicefrac{\delta}{\varepsilon}$ in the general monotone case and $\nicefrac{\delta}{\mu}\cdot\log\nicefrac{1}{\varepsilon}$ in the strongly monotone case. In contrast to existing communication-efficient algorithms that exploit similarity, \textsc{PAUS} is able to tackle non-Euclidean problems because it uses the Bregman setup.

 \newpage


\appendix

\section{Appendix with missing proofs}

\subsection{Proof of Lemma \ref{descent_lemma}}\label{app:descent_lemma}
\begin{lemma}\textbf{(Lemma \ref{descent_lemma})}
    Consider Assumption \ref{ass:delta}. Then the inequality
    \begin{align*}
  2\gamma \left[\langle F(u^k,\xi^k), u^k - z \rangle + g(u^k)-g(z)\right] \leq& 2V(z, z^k) - 2V(u^k, z^k) - 2(1+\alpha)V(z,z^{k+1}) \notag \\
  &- 2V(z^{k+1}, u^k)+2\alpha V(z,u^k) \notag+\gamma^2\delta^2\|u^k-z^k\|^2 \notag\\&+\|z^{k+1}-u^k\|^2
    \end{align*}
    holds.
\end{lemma}

\paragraph{Proof:}
\textbf{Step 1.} 
We employ the following identity for the Bregman divergence:
 \begin{align*}
     V(z, z^k) - V(z, u^{k}) &- V(u^{k}, z^k) = \langle \nabla w(u^{k}) - \nabla w (z^k), z - u^{k}) \rangle.
 \end{align*}
 Using this identity for Line \ref{line:SPP_prob} of Algorithm \ref{Alg:VI}, we get
 \begin{align*}
     \gamma(\langle F_1(u^k)+F(z^k,\xi^k)-F_1(z^k),z-u^k \rangle &+ g(z)-g(u^k))\\&+V(z,z^k)-V(z,u^k)-V(u^k,z^k)\geq0.
 \end{align*}
 We rewrite this and obtain
\begin{align}\label{eq:step1}
    \gamma \langle F_1(u^k) &+ F(z^k,\xi^k) - F_1 (z^k), u^k - z  \rangle +\gamma(g(u^k)-g(z))\leq  V(z, z^k) - V(z, u^{k}) - V(u^{k}, z^k).
\end{align}
 \textbf{Step 2.} Writing the optimality condition  for Line \ref{eq:line5Paus} of Algorithm \ref{Alg:VI},
 we get for all $z \in \Z$
 \begin{align}\label{eq:optimal_cond_row2}
  \langle  \gamma (F(u^k,\xi^k) &-F_1(u^k)-F(z^k,\xi^k)+F_1(z^k))+(1+\alpha)(\nabla w(z^{k+1}) - \nabla w (u^k)), z - z^{k+1} \rangle \geq 0.
 \end{align}
 Then we utilize the following identity for the Bregman divergence
 \begin{align*}
     V(z, u^k) - &V(z, z^{k+1}) - V(z^{k+1}, u^k) = \langle \nabla w(z^{k+1}) - \nabla w (u^k), z - z^{k+1}) \rangle.
 \end{align*}
 Using this identity
 for \eqref{eq:optimal_cond_row2} and denoting $H(u^k,z^k,\xi^k)=F(u^k,\xi^k) -F_1(u^k) 
    -F(z^k,\xi^k)+F_1(z^k)$, we obtain
\begin{align*}
   \gamma \langle H(u^k,z^k,\xi^k), z^{k+1}-z \rangle \leq& (1+\alpha)V(z, u^k) - (1+\alpha)V(z, z^{k+1}) - V(z^{k+1}, u^k).
\end{align*}
Let us add and then subtract $u^k$ in $\langle H(u^k,z^k,\xi^k), z^{k+1}-z \rangle$:
\begin{align}\label{eq:step2}
\begin{split}
   \gamma \langle H(u^k,z^k,\xi^k), z^{k+1}-u^k \rangle + \gamma \langle H(u^k,z^k,\xi^k), u^k-z \rangle \leq& (1+\alpha)V(z, u^k) \\&- (1+\alpha)V(z, z^{k+1}) \\&- V(z^{k+1}, u^k).
\end{split}
\end{align}
\textbf{Step 3.} Now we summarize \eqref{eq:step1} and \eqref{eq:step2} and get
\begin{align*}
    \gamma \langle F(u^k,\xi^k), u^k - z \rangle + \gamma(g(u^k)-g(z)) \leq& V(z, z^k) - V(u^k, z^k)-(1+\alpha)V(z,z^{k+1}) \notag\\&- V(z^{k+1}, u^k)+\alpha V(z,u^k) \notag\\&+ \gamma \langle H(u^k,z^k,\xi^k), u^k-z^{k+1}  \rangle. 
\end{align*}
 Next for this we use the Cauchy–Schwarz inequality
 \begin{align}\label{eq:cauchy-sch}
 \begin{split}
  2\gamma \left[\langle F(u^k,\xi^k), u^k - z \rangle + g(u^k)-g(z)\right] \leq& 2V(z, z^k) - 2V(u^k, z^k) - 2(1+\alpha)V(z,z^{k+1}) \\
  &- 2V(z^{k+1}, u^k)+2\alpha V(z,u^k) \\&+\gamma^2\|H(u^k,z^k,\xi^k)\|_*^2\\&+\|z^{k+1}-u^k\|^2. 
  \end{split}
 \end{align}
Assumption \ref{ass:delta} gives
\[
\|H(u^k,z^k,\xi^k)\|_*=\|F(u^k,\xi^k) - F_1(u^k)-F(z^k,\xi^k)+F_1(z^k)\|_* \leq \delta \|u^k-z^k\|.
 \]
We use this for \eqref{eq:cauchy-sch} and obtain
 \begin{align*}
  2\gamma \left[\langle F(u^k,\xi^k), u^k - z \rangle + g(u^k)-g(z)\right] \leq& 2V(z, z^k) - 2V(u^k, z^k) - 2(1+\alpha)V(z,z^{k+1}) \notag \\
  &- 2V(z^{k+1}, u^k)+2\alpha V(z,u^k) \notag+\gamma^2\delta^2\|u^k-z^k\|^2 \notag\\&+\|z^{k+1}-u^k\|^2. 
 \end{align*}
This transition completes the proof of the lemma.
 \hfill $ \square$

\subsection{Proof of Theorem \ref{th:PAUS} }\label{app:lemma_PAUS}
\begin{theorem}\textbf{(Theorem \ref{th:PAUS})}
Consider assumptions of Lemma \ref{descent_lemma} with Assumptions \ref{ass:stochastic}(a) and \ref{ass:monotone}. Then after $K $ communication rounds, \textsc{PAUS} (Algorithm \ref{Alg:VI}),  run with $\alpha=0$, stepsize $\gamma \leq\nicefrac{1}{2\delta}$ and a starting point $z^0\in \Z$,
outputs  $\tilde{u}^K$ such that
\begin{align*}
    \E\left[\Gap(\tilde{u}^K)\right]\leq\frac{D}{K\gamma}+ \frac{2\gamma}{3}\sigma^2,
\end{align*}
where $D=\sup_{z\in Z}\left\{V(z, z^0)\right\}$.
\end{theorem}
\paragraph{Proof:}
From Lemma \ref{descent_lemma} we have
 \begin{align*}
  2\gamma \left[\langle F(u^k,\xi^k), u^k - z \rangle + g(u^k)-g(z)\right] \leq& 2V(z, z^k) - 2V(u^k, z^k) - 2(1+\alpha)V(z,z^{k+1}) \notag \\
  &- 2V(z^{k+1}, u^k)+2\alpha V(z,u^k) \notag+\gamma^2\delta^2\|u^k-z^k\|^2 \notag\\&+\|z^{k+1}-u^k\|^2. 
 \end{align*}
Next we take $\alpha=0$ and utilize the monotonicity of $F(\cdot,\xi)$ (Assumption \ref{ass:monotone}):
 \begin{align*}
  2\gamma \left[\langle F(z,\xi^k), u^k - z \rangle + g(u^k)-g(z)\right] \leq& 2V(z, z^k) - 2V(u^k, z^k) - 2(1+\alpha)V(z,z^{k+1}) \notag \\
  &- 2V(z^{k+1}, u^k)+2\alpha V(z,u^k) \notag+\gamma^2\delta^2\|u^k-z^k\|^2 \notag\\&+\|z^{k+1}-u^k\|^2. 
 \end{align*}
We add and subtract $F(z)$ in the scalar product and obtain
\begin{align*}
   2\gamma \left[\langle F(z), u^k - z \rangle + g(u^k)-g(z)\right] 
  \leq& 2V(z, z^k) - 2V(u^k, z^k) -2V(z,z^{k+1}) - 2V(z^{k+1},  u^k) \notag \\
  &+\gamma^2 \delta^2\|u^k - z^k\|^2+\|z^{k+1}-u^k\|^2.
  \\
  &+ 2\gamma \langle F(z)-F(z,\xi^k), u^k - z \rangle.
\end{align*}
Note that:
\begin{align*}
    \E_{\xi}\langle F(z)-F(z,\xi^k), u^k - z\rangle = \E_{\xi}\langle F(z)-F(z,\xi^k), u^k - z^k\rangle,
\end{align*}
because $F(\cdot,\xi)$ is unbiased and $z,z^k$ are independent on $\xi^k$. Thus, we calculate the expectation and get
\begin{align*}
     2\gamma \E\langle F(z), u^k - z \rangle + g(u^k)-g(z) 
  \leq& \E2V(z, z^k) - 2V(u^k, z^k) -2V(z,z^{k+1}) \notag \\
  &- 2V(z^{k+1},  u^k) +\gamma^2 \delta^2\|u^k - z^k\|^2+\|z^{k+1}-u^k\|^2
  \\
  &+ \frac{4\gamma^2}{3}\sigma_z^2 + \frac{3}{4}\|u^k-z^k\|^2. 
\end{align*}
Here we also apply the Cauchy-Schwartz inequality to $2\gamma \langle F(z)-F(z,\xi^k), u^k - z^k \rangle$. For the Bregman divergence the inequality $V(x,y)\geq \frac{1}{2}\|x-y\|^2$ is satisfied for all $x, y$. Thus,
\begin{align*}
    -2V(u^k,z^k)\leq-\|u^k-z^k\|^2,
\end{align*}
\begin{align*}
    -2V(z^{k+1},u^k)\leq-\|z^{k+1}-u^k\|^2.
\end{align*}
Using this inequalities, we write
\begin{align*}
     2\gamma \E\langle F(z), u^k - z \rangle + g(u^k)-g(z) 
  \leq& \E2V(z, z^k) -2V(z,z^{k+1}) +\left(\gamma^2\delta^2+\frac{3}{4}-1\right)\|u^k-z^k\|^2 \\&+\left(1-1\right)\|z^{k+1}-u^k\|^2 + \frac{4\gamma^2}{3}\sigma^2.
\end{align*}
If we choose stepsize $\gamma \leq \frac{1}{2\delta}$, we obtain
\begin{align*}
     2\gamma \E\langle F(z), u^k - z \rangle + g(u^k)-g(z) 
  &\leq \E2V(z, z^k) -2V(z,z^{k+1}) + \frac{4\gamma^2}{3}\sigma^2.
\end{align*}
We summarize this for $k=0,1,2,\dots, K-1$
\[
     2\gamma\E\frac{1}{K} \sum_{k=0}^{K-1} \left[\langle F(z), u^k - z \rangle + g(u^k)-g(z) \right]  
     \leq \frac{2}{K}V(z, z^0) + \frac{4\gamma^2}{3}\sigma^2. 
\]
The statement of the theorem follows by taking the maximum:
\begin{align*}
    \E[\Gap(\tilde{u}^k)]\leq\frac{1}{K\gamma}\sup_{z\in Z}\left\{V(z, z^0)\right\} + \frac{2\gamma}{3}\sigma^2.
\end{align*}

 \hfill $ \square$

\subsection{proof of Corollary \ref{cor:monotone}}\label{proof:cor:monotone}
\begin{corollary}\textbf{(Corollary \ref{cor:monotone})}
 Consider assumptions of Theorem \ref{th:PAUS}.
Let $\tilde{u}^K$ be the output of \textsc{PAUS} (Algorithm \ref{Alg:VI}), run  with an appropriate parameters and starting points $z^0,u^0\in \Z$ 
in
 \[ \mathcal{O}\left( \frac{D\delta}{\varepsilon}+\frac{D\sigma^2}{\varepsilon^2} \right) \text{communication rounds,}
  \]  
then  $\Gap(\tilde{u}^K)\leq \varepsilon$. 
\end{corollary}

\paragraph{Proof:}
From Theorem \ref{th:PAUS} we have
\begin{align*}
    \E\left[\Gap(\tilde{u}^K)\right]\leq\frac{D}{K\gamma}+ \frac{2\gamma}{3}\sigma^2,
\end{align*}
where $D=\sup_{z\in Z}\left\{V(z, z^0)\right\}$.
Let us find $\gamma$, equating the summands in the right-hand side:
\begin{align*}
    \gamma=\frac{1}{\sigma}\sqrt{\frac{3D}{2K}}.
\end{align*}
\begin{itemize}
    \item If $\frac{1}{\sigma}\sqrt{\frac{3D}{2K}}\leq\frac{1}{2\delta}$, choose $\gamma=\frac{1}{\sigma}\sqrt{\frac{3D}{2K}}$. In this case we obtain
    \begin{align*}
         \E\left[\Gap(\tilde{u}^K)\right]\leq2\sigma\sqrt{\frac{2D}{3K}}.
    \end{align*}
    \item If $\frac{1}{\sigma}\sqrt{\frac{3D}{2K}}\geq\frac{1}{2\delta}$, choose $\gamma=\frac{1}{2\delta}$. In this case we have
    \begin{align*}
         \E\left[\Gap(\tilde{u}^K)\right]\leq \frac{2\delta D}{K} + \frac{2\sigma^2}{3}\frac{2}{2\delta}\leq \frac{2\delta D}{K} + \frac{2\sigma^2}{3}\frac{1}{\sigma}\sqrt{\frac{3D}{2K}}.
    \end{align*}
\end{itemize}
Getting rid of unnecessary constants, we get the communication complexity \[ \mathcal{O}\left( \frac{D\delta}{\varepsilon}+\frac{D\sigma^2}{\varepsilon^2} \right).
  \]  
 \hfill $ \square$

\subsection{Proof of Theorem \ref{th:compositeMP} }\label{app:irror_comp}
\begin{theorem}\textbf{(Theorem \ref{th:compositeMP})}
   Let Assumptions \ref{ass:stochastic_server}, \ref{ass:monotone}, \ref{ass:Lipsch} and \ref{ass:delta} hold. Then \textsc{SCMP} (Algorithm \ref{Alg:CMP}) with $\gamma=\frac{1}{2\delta}$ produce the sequence $\{v^t\}$ such that 
   \begin{align*}
    \E\left[V(v^*, v^{t+1})\right] \leq \left(1-\frac{\eta}{2}\right)\E\left[V(v^*,v^t)\right] + 4\eta^2\sigma_{1,*}^2.
    \end{align*}
\end{theorem}
\paragraph{Proof:}
Let us introduce 
\[
H(x,\xi) \defeq \gamma (F_1(x,\xi)+F(z^k,\xi^k) - F_1(z^k)).
\]
Then we can rewrite the iterates of \textsc{SCMP} procedure: 
\begin{align}
    v^{t+\frac{1}{2}} &\gets \arg\min\limits_{v\in \Z} \{ \eta\langle H(v^t,\xi^t), v  \rangle 
  + \eta V(v, z^k) + V(v, v^t)+\gamma g(v)\}. \label{eq:gkdkkkdkmmm} \\
  v^{t+1} &\gets \arg\min\limits_{v\in \Z} \{ \eta\langle H(v^{t+\frac{1}{2}}.\xi^t), v  \rangle 
  + \eta V(v, z^k) + V(v, v^t)+\gamma g(v)\}.\label{eq:gkdkkkdkmmm2}
\end{align}
From the optimality conditions for \eqref{eq:gkdkkkdkmmm} and \eqref{eq:gkdkkkdkmmm2} we have
\begin{align}
   \langle  \eta H(v^t,\xi^t) +& \eta (\nabla w(v^{t+\frac{1}{2}}) - \nabla w(z^k)) + \nabla w(v^{t+\frac{1}{2}}) - \nabla w(v^t) , v^{t+\frac{1}{2}} - v\rangle \notag\\&\leq \gamma(g(v)-g(v^{t+\frac{1}{2}})). \label{eq:align_gg}
\end{align}
\begin{align}
    \langle  \eta H(v^{t+\frac{1}{2}},\xi^t) +& \eta (\nabla w(v^{t+1}) - \nabla w(z^k))+ \nabla w(v^{t+1}) - \nabla w(v^t) , v^{t+1} - v\rangle \notag\\&\leq \gamma(g(v)-g(v^{t+1})). \label{eq:align_gg2}
\end{align}
Let $v^*$ be an exact solution of problems in Line \ref{line:SPP_prob} of \textsc{PAUS} for which we employ the  \textsc{SCMP} procedure.   Plugging $v = v^{t+1}$ in  \eqref{eq:align_gg} and $v = v^*$ in \eqref{eq:align_gg2} and the summarizing these, we get 
\begin{align*}
    \langle\eta H(v^t,\xi^t) +& \eta(\nabla w(v^{t+\frac{1}{2}})-\nabla w(z^k))+\nabla w(v^{t+\frac{1}{2}})-\nabla w(v^t),v^{t+\frac{1}{2}}-v^{t+1}\rangle
    \\
    +&\langle\eta H(v^{t+\frac{1}{2}},\xi^t)+\eta(\nabla w(v^{t+1})-\nabla w(z^k))+\nabla w(v^{t+1})-\nabla w(v^t),v^{t+1}-v^*\rangle
    \\
    &\leq \gamma(g(v^*)-g(v^{t+\frac{1}{2}})).
\end{align*}
Let us add and subtract $v^{t+\frac{1}{2}}$ into the second inner product and obtain
\begin{align*}
    \langle\eta H(v^{t+\frac{1}{2}},\xi^t)+&\eta(\nabla w(v^{t+1})-\nabla w(z^k)),v^{t+\frac{1}{2}}-v^*\rangle \\+& \langle \eta(H(v^t,\xi^t)-H(v^{t+\frac{1}{2}},\xi^t)),v^{t+\frac{1}{2}}-v^{t+1} \rangle 
    \\
    +&\eta\langle \nabla w(v^{t+\frac{1}{2}})-\nabla w(v^{t+1}),v^{t+\frac{1}{2}}-v^{t+1} \rangle \\+& \langle \nabla w(v^{t+\frac{1}{2}})-\nabla w(v^t), v^{t+\frac{1}{2}}-v^{t+1} \rangle 
    \\
    +& \langle \nabla w(v^{t+1})-\nabla w(v^t), v^{t+1}-v^* \rangle
    \\
    &\leq \gamma(g(v^*)-g(v^{t+\frac{1}{2}})).
\end{align*}
After rearranging the terms we obtain
\begin{align*}
    \langle\eta H(v^{t+\frac{1}{2}},\xi^t)+&\eta(\nabla w(v^{t+1})-\nabla w(z^k)),v^{t+\frac{1}{2}}-v^*\rangle
    \\
    &\leq \gamma(g(v^*)-g(v^{t+\frac{1}{2}})) \\&\quad+ \eta\langle H(v^{t+\frac{1}{2}},\xi^t)-H(v^t,\xi^t),v^{t+\frac{1}{2}}-v^{t+1} \rangle 
    \\
    &\quad+\eta\langle \nabla w(v^{t+1})-\nabla w(v^{t+\frac{1}{2}}),v^{t+\frac{1}{2}}-v^{t+1} \rangle \\&\quad+ \langle \nabla w(v^{t+\frac{1}{2}})-\nabla w(v^t), v^{t+1}-v^{t+\frac{1}{2}} \rangle 
    \\
    &\quad+ \langle \nabla w(v^{t+1})-\nabla w(v^t), v^*-v^{t+1} \rangle.
\end{align*}
Assumption \ref{ass:monotone} implies
\begin{align}\label{proof:two_sampling}
    \langle H(v^{t+\frac{1}{2}},\xi^t),v^{t+\frac{1}{2}}-v^*\rangle\geq&\langle H(v^*,\xi^t),v^{t+\frac{1}{2}}-v^*\rangle
    \notag\\=& \langle H(v^*),v^{t+\frac{1}{2}}-v^*\rangle + \langle F_1(v^*,\xi^t)-F_1(v^*),v^{t+\frac{1}{2}}-v^*\rangle.
\end{align}
Moreover, we use the optimality condition for problem in Line \ref{line:SPP_prob} of \textsc{PAUS}, and we obtain the following for all $v\in \Z$
\begin{align*}
    \langle H(v^*) +\nabla w(v^*) - \nabla w(z^k), v - v^* \rangle \geq \gamma(g(v^*)-g(v)).
\end{align*}
Plugging $v = v^{t+\frac{1}{2}}$ in this, we have
\begin{align}\label{eq:deddffdfd}
   -\langle H(v^*)+\nabla w(v^*),v^{t+\frac{1}{2}}-v^* \rangle + \langle  \nabla w(z^k),v^{t+\frac{1}{2}}-v^*\rangle\leq \gamma(g(v^{t+\frac{1}{2}})-g(v^*)).  
\end{align}
  Thus, summarizing \eqref{eq:deddffdfd} and the origin inequality, we get 
  \begin{align}\label{eq:ghhdddfhj}
     \eta\langle F_1(v^*,\xi^t)-F_1(v^*),v^{t+\frac{1}{2}}-v^* \rangle +& \eta \langle \nabla w(v^{t+1}) - \nabla w(v^*), v^{t+\frac{1}{2}} - v^*\rangle \notag \\
      &\leq \eta \langle H(v^{t+\frac{1}{2}},\xi^t) - H(v^t,\xi^t), v^{t+\frac{1}{2}} - v^{t+1}\rangle\notag  \\&\quad+\eta \langle \nabla w(v^{t+1})- \nabla w(v^{t+\frac{1}{2}}), v^{t+\frac{1}{2}} - v^{t+1}  \rangle \notag \\
        &\quad+ \langle \nabla w(v^{t+\frac{1}{2}}) - \nabla w(v^t) , v^{t+1} - v^{t+\frac{1}{2}} \rangle\notag \\&\quad+ \langle \nabla w(v^{t+1}) - \nabla w(v^t), v^* - v^{t+1}\rangle.
  \end{align}

Note that $\E_{\xi^t}[\langle F_1(v^*,\xi^t)-F_1(v^*),v^{t+\frac{1}{2}}-v^* \rangle]\neq0$, because $v^{t+\frac{1}{2}}$ depends on $\xi^t$. Assumption \ref{ass:stochastic_server} allows to replace $v^*$ by $v^t$ in the first summand of the left-hand side while taking expectation:
\begin{align*}
    \E_{\xi^t}\langle F_1(v^*,\xi^t)-F_1(v^*),v^{t+\frac{1}{2}}-v^* \rangle =& \E_{\xi^t}\langle F_1(v^*,\xi^t)-F_1(v^*),(v^{t+\frac{1}{2}}-v^t)-(v^*-v^t) \rangle.
\end{align*}
Both $v^*$ and $v^t$ are independent on $\xi^t$, Thus, we obtain
\begin{align*}
    \E_{\xi^t}\langle F_1(v^*,\xi^t)-F_1(v^*),v^{t+\frac{1}{2}}-v^* \rangle = \E_{\xi^t}\langle F_1(v^*,\xi^t)-F_1(v^*),v^{t+\frac{1}{2}}-v^t \rangle.
\end{align*}
Hence, using this, we rewrite \eqref{eq:ghhdddfhj} as follows
\begin{align}\label{eq:fjpsjjnnn}
   \eta\E \langle  \nabla w(v^{t+1}) - \nabla w(v^*), v^{t+\frac{1}{2}} -v^*\rangle\leq&  \E\eta\langle F_1(v^*,\xi^t)-F_1(v^*),v^t-v^{t+\frac{1}{2}} \rangle \notag\\&+ \eta \langle H(v^{t+\frac{1}{2}}, \xi^t) - H(v^t,\xi^t), v^{t+\frac{1}{2}} - v^{t+1}  \rangle \notag 
   \\
   &+ \eta \langle \nabla w(v^{t+1}) - \nabla w(v^{t+\frac{1}{2}}), v^{t+\frac{1}{2}} - v^{t+1} \rangle \notag\\&+ \langle \nabla w(v^{t+\frac{1}{2}}) - \nabla w(v^t), v^{t+1}  - v^{t+\frac{1}{2}}\rangle \notag
   \\
   &+ \langle \nabla w(v^{t+1}) - \nabla w(v^{t}), v^* - v^{t+1} \rangle.
\end{align}
From the definition of the Bregman divergence, we have 
\begin{align}
    - V(v^{t+1}, v^{t+\frac{1}{2}}) - V(v^{t+\frac{1}{2}}, v^{t+1}) &= \langle \nabla w(v^{t+1}) - \nabla w(v^{t+\frac{1}{2}}), v^{t+\frac{1}{2}} - v^{t+1} \rangle. \label{eq:div_breg1} \\
      V(v^*, v^{t+1}) + V(v^{t+1}, v^{t})  - V(v^*,v^{t})&= \langle \nabla w(v^t) - \nabla w(v^{t+1}), v^* - v^{t+1}\rangle.
    \label{eq:div_breg2} \\
    V(v^{t+1}, v^{t+\frac{1}{2}}) + V(v^{t+\frac{1}{2}}, v^{t})  - V(v^{t+1},v^{t})&= \langle \nabla w(v^{t}) - \nabla w(v^{t+\frac{1}{2}}), v^{t+1} - v^{t+\frac{1}{2}} \rangle. \label{eq:div_breg3} \\
        V(v^{t+\frac{1}{2}}, v^*) + V(v^*, v^{t+1})  - V(v^{t+\frac{1}{2}},v^{t+1})&= \langle \nabla w(v^{t+1}) - \nabla w(v^*), v^{t+\frac{1}{2}} - v^* \rangle.\label{eq:div_breg4} 
\end{align}
Plugging \eqref{eq:div_breg1}, \eqref{eq:div_breg2}, \eqref{eq:div_breg3} and \eqref{eq:div_breg4} in \eqref{eq:fjpsjjnnn}, we obtain 
\begin{align*}
   \E\eta V(v^{t+\frac{1}{2}}, v^*) + \eta V(v^*, v^{t+1})  - \eta V(v^{t+\frac{1}{2}},v^{t+1}) 
   \leq& \E\eta\langle F_1(v^*,\xi^t)-F_1(v^*),v^t-v^{t+\frac{1}{2}} \rangle
   \\
   &+\eta \langle H(v^{t+\frac{1}{2}},\xi^t) - H(v^t,\xi^t), v^{t+\frac{1}{2}} - v^{t+1} \rangle \notag \\
   &- \eta V(v^{t+1}, v^{t+\frac{1}{2}}) - \eta V(v^{t+\frac{1}{2}}, v^{t+1}) \notag \\
   &   - V(v^{t+1}, v^{t+\frac{1}{2}}) - V(v^{t+\frac{1}{2}}, v^t) \notag\\&+ V(v^{t+1}, v^t)
    -V(v^*, v^{t+1}) - V(v^{t+1}, v^t) \notag\\&+V(v^*, v^t).
\end{align*}
Rearranging the terms and using $V(v^{t+\frac{1}{2}}, v^*) \geq 0$, we get
\begin{align}\label{eq:ggkll;d}
    0
    \leq& \E\eta\langle F_1(v^*,\xi^t)-F_1(v^*),v^t-v^{t+\frac{1}{2}} \rangle + \eta \langle H(v^{t+\frac{1}{2}},\xi^t) - H(v^t,\xi^t), v^{t+\frac{1}{2}} - v^{t+1} \rangle \notag\\
    &- (1+\eta) V(v^{t+1}, v^{t+\frac{1}{2}})- V(v^{t+\frac{1}{2}}, v^t) - (1+\eta)V(v^*,v^{t+1}) +V(v^*,v^t).
\end{align}
By using the Cauchy–Schwarz inequality with some constants $C_1$ and $C_2$, we obtain
\begin{align*}
   (1+\eta)\E V(v^*,v^{t+1})\leq&\E V(v^*,v^t)-(1+\eta)V(v^{t+1},v^{t+\frac{1}{2}})-V(v^{t+\frac{1}{2}},v^t) \\&+ C_1\sigma_*^2+\frac{1}{C_2}\|v^{t+\frac{1}{2}}-v^{t+1}\|^2+\left(\frac{\eta^2}{C_1}+C_2\eta^2\gamma^2L_{F_1}^2\right)\|v^{t+\frac{1}{2}}-v^t\|^2.
\end{align*}

Next we use the fact that $V(x,y) \geq \frac{1}{2}\|x-y\|^2$ for all $x, y \in \R^d$. Let us choose
\begin{equation}\label{eq:etaC}
    \gamma=\frac{1}{2\delta}, C_1 = 4\eta^2, C_2=2, \eta \leq \frac{1}{3\gamma L_{F_1}}.
\end{equation}
$\frac{1}{1+\eta} \leq 1-  \frac{x}{2}$, since $\eta\leq\frac{1}{3\gamma L_{F_1}}=\frac{2\delta}{3 L_{F_1}}<1$. Thus, we have
\begin{equation*}
    \E[V(v^*, v^{t+1})] \leq \left(1-\frac{\eta}{2}\right)\E[V(v^*,v^t)] + 4\eta^2\sigma_*^2.
\end{equation*}

\subsection{Proof of Corollary \ref{cor:scmp}}\label{proof:cor:scmp}
\begin{corollary}\textbf{(Corollary \ref{cor:scmp})}
   Consider assumptions of Theorem \ref{th:compositeMP}. Let  $v^*$ be a solution  of the subproblem in Line \ref{line:SPP_prob}  of Algorithm \ref{Alg:CMP} and
   let $ v^{T}$  be the output of  \textsc{SCMP} procedure. Consider stepsize $\gamma=\nicefrac{1}{2\delta}$ and starting point $v^0$. Then Algorithm \ref{Alg:CMP} with appropriate choice of $\eta$ needs
   $$
   \mathcal{O}\left(\frac{L_{F_1}}{\delta}\log\frac{V(v^*,v^0)}{\varepsilon} + \frac{\sigma_{1,*}^2}{\varepsilon}\right)
   \text{ iterations} 
   $$ 
   to achieve $V(v^*, v^T)\leq \varepsilon$.
\end{corollary}
\paragraph{Proof:}
Denote $a=\frac{1}{2}$ and $c=4\sigma_*^2$. Using \citep{stich2019unified}, we obtain:
\begin{align*}
    \E V(v^*,v^T)\leq \frac{2V(v^*,v^0)}{\eta}\exp\left\{-\frac{1}{2}\eta(T+1)\right\}+8\eta\sigma_*^2.
\end{align*}
$\bullet$ If $\frac{\delta}{3L_{F_1}}\geq \nicefrac{\ln\left(\max\left\{2,\frac{V^0T^2}{16\sigma_*^2}\right\}\right)}{T}$, then choose $\eta=\nicefrac{2\ln\left(\max\left\{2,\frac{V^0T^2}{16\sigma_*^2}\right\}\right)}{T}$ and obtain that the right side is $\mathcal{O}\left(\frac{8\sigma_*^2}{T}\right)$.\\
$\bullet$ Otherwise, choose $\eta=\frac{2\delta}{3L_{F_1}}$ and obtain $\mathcal{O}\left( \frac{3L_{F_1}V^0}{2\delta}\exp\{-\frac{\delta T}{3L_{F_1}}\}  + \frac{8\sigma_*^2}{T}\right)$.

\subsection{Proof of Theorem 
\ref{th:strongly_monotone}}\label{app:strongly}

\begin{theorem}\textbf{(Theorem \ref{th:strongly_monotone})}
Consider assumptions of Lemma \ref{descent_lemma} with Assumption \ref{ass:stochastic}(b) and Assumption \ref{strong_monotonicity}.
Consider $\alpha=\nicefrac{\gamma\mu}{2}$, $\gamma \leq \nicefrac{1}{2\delta}$ and a starting point $z^0\in \Z$. Then the inequality
\begin{align*}
    \E\left[V(z^*,z^{k+1})\right]\leq\left(1-\frac{\gamma\mu}{4}\right)\E\left[V(z^*,z^k)\right]+\frac{2\gamma^2}{3}\sigma_*^2
\end{align*}
holds.
\end{theorem}

\paragraph{Proof:}

Let us start with Lemma \ref{descent_lemma}:
\begin{align}\label{aoa}
    2\gamma \left[\langle F(u^k,\xi^k), u^k - z \rangle + g(u^k)-g(z)\right] \leq& 2V(z, z^k) - 2V(u^k, z^k) - 2(1+\alpha)V(z,z^{k+1}) \notag \\
    &- 2V(z^{k+1}, u^k)+2\alpha V(z,u^k) \notag+\gamma^2\delta^2\|u^k-z^k\|^2 \notag\\&+\|z^{k+1}-u^k\|^2. 
\end{align}
Write down the optimality condition for problem (\ref{eq:VI}):
\begin{align}\label{oao}
    \langle F(z^*),z-z^* \rangle \geq g(z^*)-g(z),\quad\forall z\in Z
\end{align}
Take $z=z^*$ in (\ref{aoa}) and $z=u^k$ in (\ref{oao}). By summing these two expressions and then adding and subtracting $\langle F(z^*,\xi^k), u^k-z^* \rangle$, we obtain the following:
\begin{align*}
   2\gamma\langle F(u^k,\xi^k)-F(z^*,\xi^k),u^k-z^* \rangle \leq& 2V(z^*, z^k) - 2V(u^k, z^k) - 2(1+\alpha)V(z^*,z^{k+1}) \notag \\
    &- 2V(z^{k+1}, u^k)+2\alpha V(z^*,u^k) \notag+\gamma^2\delta^2\|u^k-z^k\|^2 \notag\\&+\|z^{k+1}-u^k\|^2 + 2\gamma\langle F(z^*)-F(z^*,\xi^k),u^k-z^* \rangle. 
\end{align*}
Again we use the trick of replacing the point independent of $\xi^k$ by an arbitrary point independent of $\xi^k$, under the expectation due to Assumption \ref{ass:stochastic}:
\begin{align*}
    \E_{\xi^k}\langle F(z^*)-F(z^*,\xi^k),u^k-z^* \rangle = \E_{\xi^k}\langle F(z^*)-F(z^*,\xi^k),u^k-z^k \rangle.
\end{align*}
Let us apply Young's inequality to $\E_{\xi^k}\langle F(z^*)-F(z^*,\xi^k),u^k-z^k \rangle$ and Assumption \ref{strong_monotonicity} to $\langle F(u^k,\xi^k)-F(z^*,\xi^k),u^k-z^* \rangle$:
\begin{align*}
   \E\gamma\mu V(z^*,u^k) &\leq \E2V(z^*,z^k) + 2\alpha V(z^*,u^k) - 2(1+\alpha)V(z^*,z^{k+1}) + \frac{4\gamma^2}{3}\sigma_*^2.
\end{align*}
Since $\alpha=\nicefrac{\gamma\mu}{2}$, $V(z^*,u^k)$ is reduced. Note that $\nicefrac{\gamma\mu}{2}\leq\nicefrac{\mu}{4\delta}<1$. Thus, we obtain
\begin{align*}
    \E\left[V(z^*,z^{k+1})\right]\leq\left(1-\frac{\gamma\mu}{4}\right)\E\left[V(z^*,z^k)\right]+\frac{2\gamma^2}{3}\sigma_*^2
\end{align*}

\hfill $ \square$

\section{Closed forms for monotone VIs}\label{closed_forms_app}
For simplicity of presentation, we consider a non-stochastic version of \textsc{PAUS}.
\paragraph{Convex minimization.} For convex minimization problem, operator $F(z)= \nabla f(z)$, $Q= \nabla f(z) - f_1(z) =\vcentcolon \nabla q(z)$. At each iteration of of \textsc{PAUS} server forms the  gradient 
 by averaging local gradients calculated by all machines 
and then  computes the next iterates $z^{k+1}$ and $u^k$ as follows 
\begin{align}
       &u^k = \arg\min\limits_{z \in \Z} \{ \gamma f_1(z)  +\gamma \la \nabla q(z^k), z \ra + V(z, z^k)+\gamma g(z) \}, \label{eq:convmin_paus} \\
       & z^{k+1} =\arg\min\limits_{z \in \Z} \{ \gamma \la \nabla q(u^k)  -\nabla q(z^k), z \ra + V(z,u^k) \}.  \label{eq:convmin_paus2} 
\end{align}
Then the server broadcasts $z^{k+1}$ and $u^k$ to all other devices. 

In the entropy setup when $\Z \equiv \Delta$ the inner problem encountered in \eqref{eq:convmin_paus2}  has a closed-form solution  known as entropic mirror descent   \cite{nemirovski2004prox}:
 \begin{align}
     z^{k+1}&=  \frac{u^k \odot e^{-\gamma \left(\nabla q(u^k)  -\nabla q(z^k)\right) }}{  \boldsymbol{1}^\top \left( u^k \odot e^{-\gamma  \left(\nabla q(u^k)  -\nabla q(z^k)\right)}\right)}, \label{eq:entropic1} 
 \end{align}
where $\boldsymbol{1}$ is the vector of ones, exp is applied element-wise for vectors and symbols $\odot$ and $/$ stand for the element-wise product and division respectively.

\paragraph{SPPs.}
For SPPs,  operator $F(z) = [\nabla_x f(x,y),~ - \nabla_y f(x,y)]$ and $F_1(z) = [\nabla_x f_1(x,y),~ - \nabla_y f_1(x,y)]$, $G(z) \defeq F(z) - F_1(z)$. with $z \defeq (x,y) \in \X \times Y =\vcentcolon \Z$. Then the server computes
\begin{align}
   &u^k  = \arg\min\limits_{x \in \X} \max\limits_{y \in \Y} \{ \gamma f_1(x,y)  +\gamma \left\la Q(z^k), z \right\ra + V(z, z^k)+\gamma g(z) \}, \label{eq:spps_paus}\\
      & z^{k+1} =\arg\min\limits_{z \in \Z} \{ \gamma \la Q(u^k)-Q(z^k), z \ra+ V(z,u^k) \}. \label{eq:spps_paus2}
\end{align}
Similarly to convex minimization problem, in the entropy setup ($\X \equiv \Delta$ and $\Y \equiv \Delta$) the inner problem from  \ref{eq:spps_paus2}  has a closed-form solution.
\paragraph{Closed-form solutions for subproblems in \textsc{Composite MP} procedure.} 

Next we comment 
on the existence of closed-form solutions for steps \eqref{eq:MP_line5} and \eqref{eq:MP_line6} of  the \textsc{SCMP} procedure in the Entropy setup.
Particularly for convex minimization problem \eqref{eq:SPP_empir} with  $\Z \equiv \Delta$ and $g(v)\equiv0$,  \textsc{SCMP} can be rewritten as follows:
\begin{align*}
    v^{t+\frac{1}{2}} &= \frac{(z^k)^{\frac{\eta}{\eta+1}} \odot  \left(v^t \right)^{\frac{1}{1+\eta}} \odot  e^{ - \frac{\gamma \eta }{1+\eta}h(v^t) }}
    {
    \boldsymbol{1}^\top  \left( (z^k)^{\frac{\eta}{\eta+1}} \odot \left(v^t \right)^{\frac{1}{1+\eta}} \odot e^{ - \frac{\gamma \eta }{1+\eta} h(v^t)} \right)
    },  \\
 v^{t+1}  &= \frac{(z^k)^{\frac{\eta}{\eta+1}} \odot  \left(v^{t} \right)^{\frac{1}{1+\eta}} \odot  e^{ - \frac{\gamma \eta }{1+\eta}h\left(v^{t+\frac{1}{2}}\right) }}
    {
    \boldsymbol{1}^\top  \left( (z^k)^{\frac{\eta}{\eta+1}} \odot \left(v^{t}\right)^{\frac{1}{1+\eta}} \odot e^{ - \frac{\gamma \eta }{1+\eta} h\left(v^{t+\frac{1}{2}}\right)} \right)
    }, 
\end{align*}
where   $h(v) \defeq  \nabla f_1(v)+\nabla f(z^k)  -\nabla f_1(z^k) $ and  $\boldsymbol{1}$ is the vector of ones, exp is applied element-wise for vectors and symbols $\odot$ and $/$ stand for the element-wise product and division respectively. 
Similarly  closed-form solutions can be obtained for SPP with $\X \equiv \Delta$ and $\Y \equiv \Delta$.

\section{Experiments details}\label{sec:additional_ex}
We consider a two-player matrix game
\begin{equation}\label{eq:prob_matrix_game}
  \min_{x \in \Delta} \max_{y \in \Delta} \left[x^\top \bar A  y \defeq \frac{1}{m} \sum_{i=1}^m x^\top  A_{i} y\right],
\end{equation}
where $A_{1}, \dots, A_m$ are 
i.i.d samples  of stochastic matrix  $A_{\xi}$ of size $d\times d$,  $m = 10^4$. Local datasets are of size  $n =2\cdot 10^3$, the server holds also $n$ matrices.  

Next we comment on theoretical bounds for parameters $L$  and $\delta$ for this problem since $\nicefrac{1}{L}$ is the stepsize for Mirror Prox \cite{rogozin2021decentralized}, and $\nicefrac{1}{2\delta}$ is the stepsize for \textsc{PAUS} and the Eucidean algorithm \cite{kovalev2022optimal}. 
\paragraph{Lipschitz constant $L$.}
For SPP \eqref{eq:prob_matrix_game}, Assumption \ref{ass:Lipsch} is equivalent to the notion  of smoothness. 
\begin{definition}[$L$-smoothness]\label{def:smooth}
$f(x,y)$ is $ (L_{xx},L_{xy}, L_{yx}, L_{yy})$-smooth if for any $x, x' \in \X$ and $y,y' \in \Y$, 
\begin{align*}
     \|\nabla_x f(x,y) - \nabla_x f(x',y)\|_{*}
     &\leq L_{xx}\|x-x' \|,\\
       \|\nabla_x f(x,y) - \nabla_x f(x,y')\|_{*}
       &\leq L_{xy} \|y-y' \|,\\
         \|\nabla_y f(x,y) - \nabla_y f(x,y')\|_{*}
         &\leq L_{yy}\|y-y' \| ,\\
           \|\nabla_y f(x,y) - \nabla_y f(x',y)\|_{*}
           &\leq L_{yx} \|x-x' \|.
\end{align*}
\end{definition}
Then we define $L = \max\{ L_{xx}, L_{xy}, L_{yx}, L_{yy}\}$ and seek to estimate $L$.
We equip both $\X \defeq \Delta$ and  $\Y \defeq \Delta$ with the $\ell_1$-norm. The corresponding dual norm is the $\ell_\infty$-norm.  By the Definition  \ref{def:smooth}
\begin{align}\label{eq:ggddd}
    \| \bar A(y - y') \|_\infty\leq \|\sum_{i=1}^d A^{(i)}(y_i - y'_i) \|_\infty \leq \|\bar A\|_{\max}\|y - y'\|_1,
\end{align}
where we used $A^{(i)}$ for the $i$-th column of $A$, and $\|A\|_{\max}$ for the maximal entry of $A$ (in absolute value).
Thus, $L_{xy} =  L_{yx} = \|\bar A\|_{\max}$, and $L_{xx} = L_{yy} = 0$. Hence, $L = \|\bar A\|_{\max}$.
\paragraph{$\delta$-similarity.} Now we seek to estimate $\delta$.
We use Assumption \ref{ass:delta} particularly for ${F(z)= [ \nabla_x f(x,y), ~ -\nabla_y f(x,y) ]}$
\begin{align}\label{eq:gghkfllf}
    & \left\| \begin{bmatrix}
        &\nabla_x f_1(x,y)\\
        &- \nabla_y f_{1}(x,y)
    \end{bmatrix} - \begin{bmatrix}
        &\nabla_x f(x,y)\\
        &- \nabla_y f(x,y)
    \end{bmatrix} - \begin{bmatrix}
        &\nabla_{x'} f_1(x',y')\\
        &- \nabla_{y'} f_{1}(x',y')
    \end{bmatrix} + \begin{bmatrix}
        &\nabla_{x'} f(x',y')\\
        &- \nabla_{y'} f(x',y')
    \end{bmatrix} \right\|_\infty \notag\\&\leq \delta \left\| \begin{bmatrix}
        &x\\
        &y
    \end{bmatrix} - \begin{bmatrix}
        &x'\\
        &y'
    \end{bmatrix}\right\|_1.
\end{align}
 where
global objective is
\[f(x, y ) =  x^\top  \bar A y  \defeq \frac{1}{m} \sum_{i=1}^m x^\top  A_{i} y,\]
and local (stored on the server) objective is
\[f_1(x, y ) =  x^\top  \bar A^{N_1} y  \defeq\frac{1}{N_1} \sum_{\ell=1}^{N_1} x^\top  A_{\ell} y.\] 
Using this, we can rewrite \eqref{eq:gghkfllf} as follows
\begin{align}
      & \left\| \begin{bmatrix}
        & (\bar A - \bar A^{N_1})x\\
        &- (\bar A - \bar A^{N_1})^\top y
    \end{bmatrix} - \begin{bmatrix}
        &(\bar A - \bar A^{N_1})x'\\
        &- (\bar A - \bar A^{N_1})^\top y'
    \end{bmatrix} 
    \right\|_\infty \leq \delta \left\| \begin{bmatrix}
        &x - x'\\
        &y-y'
    \end{bmatrix} \right\|_1.
\end{align}
Thus,
\begin{align}\label{eq:lflllffl}
     & \left\| \begin{bmatrix}
        & (\bar A - \bar A^{N_1})(x -x')\\
        &- (\bar A - \bar A^{N_1})^\top (y-y')
    \end{bmatrix} 
    \right\|_\infty \leq \delta \left\| \begin{bmatrix}
       &x - x'\\
        &y-y'
    \end{bmatrix}\right\|_1.
\end{align}
Let us define $z = (x,y)$ and $z' = (x',y')$. Then we rewrite \eqref{eq:lflllffl} as follows
\begin{align*}
       & \left\| 
        \mathbf A (z -z')
    \right\|_\infty \leq \delta \left\| 
      z -z'
  \right\|_1,
\end{align*}
where 
\begin{align*}
    \mathbf{A} \defeq \begin{pmatrix}
        & (\bar A - \bar A^{N_1}) & \boldsymbol 0_{d\times d }\\
        &\boldsymbol  0_{d\times d } & - (\bar A - \bar A^{N_1})^\top
    \end{pmatrix}.
\end{align*}
Here $\boldsymbol  0_{d\times d }$ is the zero matrix of size $d\times d$.
By the same arguments as in \eqref{eq:ggddd} we conclude that $\delta = \|\mathbf A\|_{\max}$.

\end{document}